\documentclass[10pt]{article}

\usepackage[draft]{graphics}
\usepackage{graphicx}          	 
\usepackage{bm}                 
\usepackage{amsmath}             
\usepackage{amssymb}
\usepackage{amsfonts}             
\usepackage{verbatim}           
\usepackage{amsthm}             
 	
\usepackage{epsfig}
\usepackage{graphicx}
	
\usepackage{mathtools}
\usepackage{relsize}

\usepackage{makeidx}

\theoremstyle{plain}             

\newtheorem{theorem}{Theorem}[section]

\theoremstyle{definition}
\newtheorem{example}[theorem]{Example}

\newtheorem{remark}[theorem]{Remark}

\makeatletter
\@addtoreset{equation}{chapter}

\makeatother

\numberwithin{equation}{section}

\def\eqref#1{(\ref{#1})}
\def\dsp{\displaystyle}
\def\Frac#1#2{\frac
{
 {\raise.6ex
 \hbox{$\displaystyle#1$}}
}
{
 {\lower.6ex
 \hbox{$\displaystyle#2$}}
 }
}

\def\bigOxe{\sqcup \kern-2.3mm \sqcap}

\def\eoexample{{\unskip\nobreak\hfil\penalty50	
\hskip2em\hbox{}\nobreak\hfil$\diamondsuit$
\parfillskip=0pt\finalhyphendemerits=0\medbreak}}

\def\eoremark{{\unskip\nobreak\hfil\penalty50	
\hskip2em\hbox{}\nobreak\hfil$\triangle$
\parfillskip=0pt\finalhyphendemerits=0\medbreak}}

\def\eps{\varepsilon}

\makeatother

\def\dsp{\displaystyle}
\def\Frac#1#2{\frac
{
 {\raise.6ex
 \hbox{$\displaystyle#1$}}
}
{
 {\lower.6ex
 \hbox{$\displaystyle#2$}}
 }
}

\input epsf

\def\CHF#1#2#3{
{}_1F_1\left(
\begin{array}{c}
\begin{array}{cc} \hskip-10pt#1 \end{array}\\
\begin{array}{c}  \hskip-10pt#2 \end{array}
\end{array}
\hskip-8pt;\,#3
\right)}

\def\CHFs#1#2#3{
{}_1F_1\left({a};{c};{z}\right)
}

\def\binom#1#2{
\renewcommand{\arraystretch}{0.9}
\left(
\begin{array}{c}
\begin{array}{c}\hskip-10pt#1\end{array}\\
\begin{array}{c}\hskip-10pt#2\end{array}
\end{array}
\hskip-10pt
\renewcommand{\arraystretch}{1.0}
\right)}

\def\erfc{{\rm erfc}}

\def\intp{\int_0^\infty}
\def\intr{\int_{-\infty}^\infty}
\def\bigO{{\cal O}}
\def\calC{{{\cal C}}}

\def\calL{{{\cal L}}}

\def\tfrac#1#2{{{\lower.6ex
\hbox{$\scriptstyle#1$}}\over 
{\raise.7ex
\hbox{$\scriptstyle#2$}}}}

\def\CC{\mathbb C}
\def\RR{\mathbb R}

\def\erf{{{\rm erf}}}
\def\erfc{{{\rm erfc}}}
\def\Ai{{{\rm Ai}}}
\def\Bi{{{\rm Bi}}}
\def\erfc{{\rm erfc}}

\def\intp{\int_0^\infty}
\def\intr{\int_{-\infty}^\infty}

\def\calL{{\cal L}}

\def\calD{{\cal D}}
\def\phase{{\rm ph}}

\def\tfrac#1#2{{{\lower.6ex
\hbox{$\scriptstyle#1$}}\over 
{\raise.7ex
\hbox{$\scriptstyle#2$}}}}

\def\sn{\sum_{n=0}^\infty\,}
\def\sk{\sum_{k=0}^\infty\,}

\def\insil#1{}

\begin{document}
 \title{
Uniform Asymptotic Methods for Integrals
}

\author{
    Nico M. Temme\\
    IAA, 1391 VD 18, Abcoude, The Netherlands\footnote{Former address: Centrum Wiskunde \& Informatica (CWI), 
        Science Park 123, 1098 XG Amsterdam,  The Netherlands}\\
    {\ }\\
     { \small e-mail: {\tt Nico.Temme@cwi.nl}}
}          

\maketitle
\begin{abstract}\noindent
We give an overview of basic methods that can be used for obtaining asymptotic expansions of integrals: Watson's lemma, Laplace's method, the saddle point method, and the method of stationary phase. Certain developments in the field of asymptotic analysis will be compared with De Bruijn's book {\em Asymptotic Methods in Analysis}. The classical methods can be modified for obtaining expansions that hold uniformly with respect to additional parameters. We give an overview of examples in which special functions, such as the complementary error function, Airy functions, and Bessel functions, are used as approximations in uniform asymptotic expansions.
\end{abstract}

\vskip 0.8cm \noindent
{\small
2000 Mathematics Subject Classification:
41A60, 30E15, 33BXX, 33CXX.
\par\noindent
Keywords \& Phrases:
asymptotic analysis,
Watson's lemma,
Laplace's method,
method of stationary phase,
saddle point method,
uniform asymptotic expansions,
special functions.
}

\def\cprime{$'$}

\section{Introduction}\label{sec:intro}

Large parameter problems occur in all branches of pure and applied mathematics, in physics and engineering, in statistics and probability theory. On many occasions the problems are presented in the form of integrals or differential equations, or both, but we also encounter finite sums,  infinite series, difference equations, and implicit algebraic equations. 

Asymptotic methods for handling these problems have a long history with many prominent contributors. In 1863, Riemann used the method of steepest descent
for hypergeometric functions, and in 1909 Debye \cite{Debye:1909:NFZ} used this method to obtain approximations
for Bessel functions. In other unpublished notes, Riemann also gave  the first steps for approximating the zeta function and  in 1932 Siegel used this method to derive the Riemann-Siegel formula for the Riemann zeta function. The method of stationary phase was essential in Kelvin's work to describe the wave pattern behind a moving ship.

In the period 1940--1950, a systematic study of asymptotic methods
started in the Netherlands, driven by J.~G.~ van der Corput, who considered these
methods important when studying problems of number theory in his earlier years in
Groningen. Van der Corput was one of the founders of the Mathematisch Centrum in Amsterdam, and in 1946 he became the first director. He organized working groups and colloquia on asymptotic methods and he had much influence on workers in this area. During 1950--1960 he wrote many papers on asymptotic methods and during his years in the USA he published lecture notes and technical reports. 

Nowadays his work on asymptotic analysis for integrals is still of interest, in particular his work on the method of stationary phase. He introduced for this topic the so-called neutralizer in order to handle integrals with several points that contribute to the asymptotic behavior of the complete integral. In later sections we give more details.

Although Van der Corput wanted\footnote{He concluded a Rouse Ball lecture at Cambridge (England) in 1948 with a request \cite{Corput:1948:MCP} to all workers in asymptotics to send information to the MC, and a study group would include  these contributions in a general survey of the whole field.} to publish a general compendium to
asymptotic methods, the MC lecture notes and his many published papers, notes and
reports were never combined into the standard work he would have liked.
Actually, there were no books on asymptotic methods before 1956. In certain books and published papers these methods were considered in great detail. 

For example, in Watson's book on Bessel functions \cite{Watson:1944:TTB} (first edition in 1922) and in Szeg{\H{o}}'s book on orthogonal polynomials \cite{Szego:1975:OP} (first edition 1939) many results on asymptotic expansions can be found. Erd\'elyi's book {\em Asymptotic Expansions} \cite{Erdelyi:1956:ASE} came in 1956, and the first edition of De Bruijn's book {\em Asymptotic Methods in Analysis} \cite{Bruijn:1981:AMA} was published in 1958. This remarkable book was based on a course of lectures at the MC in 1954/1955 and at the Technical University in Eindhoven in 1956/1957.

De Bruijn's book was well received. It was so interesting because of the many special  examples and topics, and because the book was sparklingly written. The sometimes entertaining style made Erd\'elyi remark in his comments in the Mathematical Reviews that the book was somewhat conversational. In addition, he missed the work of Van der Corput on the method of stationary phase. These minor criticisms notwithstanding, Erd\'elyi  was, in his review, very positive about this book .

De Bruijn motivated his selection of topics and way of treating them by observing that asymptotic methods are very flexible, and in such cases it is not possible to formulate a single theorem that covers all applications. It seems that he liked  the saddle point method very much. He spent three chapters on explaining the method and provided difficult examples. For example, he considered
the sum
\begin{equation}\label{eq:intro01}
S(s,n)=\sum_{k=0}^{2n}(-1)^{k+n}{2n \choose k}^s
\end{equation}  
for large values of $n$ with $s$ a real number. Furthermore, he discussed the generalization of Euler's gamma function in the form
\begin{equation}\label{eq:intro02}
G(s)=\int_0^\infty e^{-P(u)}u^{s-1}\,du, \quad P(u)=u^N+a_{N-1}u^{N-1}+\ldots +a_0
\end{equation}  
with $\Re s\to\infty$. 

For many readers De Bruijn's treatment of the saddle
point method serves as the best introduction to this topic and many publications refer
to this part of the book. Other topics are iterative methods, implicit functions, and differential equations, all with apparently simple examples, but also with sometimes complicated forms of asymptotic analysis.

In this paper we give a short description of the classical methods that were available for one-dimensional integrals when De Bruijn's book appeared in 1958: integrating by parts, the method of stationary phase,  the saddle point method, and the related method of steepest descent. Shortly after this period, more powerful methods were developed in which extra parameters (and not only the large parameter) were taken into account. In these so-called uniform methods the asymptotic approximations usually contain certain special functions, such as error functions, Airy functions, Bessel functions, and parabolic cylinder functions. Some simple results were already available in the literature. For example, Watson \cite[\S8.43]{Watson:1944:TTB} mentions results from 1910 of Nicholson  in which Airy functions have been used.

In a later section we give an overview of integrals with large and additional parameters for which we need uniform expansions. 
For a few cases we show details of
uniform expansions of some well-known special functions, for which these
uniform expansions are important

\subsection{Notation and symbols used in asymptotic estimates }\label{sec:symbols}
For information on the special functions we use in this paper, we refer to \cite{Olver:2010:HMF}, the {\em NIST Handbook of Mathematical Functions}.\footnote{A web version is available at http://dlmf.nist.gov/}

We  use Pochhammer's symbol $(\lambda)_n$, or shifted factorial, defined by $(\lambda)_0=1$ and
\begin{equation}\label{eq:not01}
 (\lambda)_n=\frac{\Gamma(\lambda+n)}{\Gamma(\lambda)}=\lambda(\lambda+1)\cdots(\lambda+n-1),\quad n\ge1.
\end{equation}

In asymptotic estimates we use the big $O$-symbol, denoted by $\bigO$, and the little $o$-symbol.

For estimating a function $f$ with respect to $g$, both functions defined in domain $\calD\in\CC$, we assume that  $g(z)\ne0, z\in\calD$ and $z_0$ is a limit point of $\calD$. Possibly $g(z)\to0$ as $z\to z_0$. We use the $\bigO$-symbol in the form
\begin{equation}\label{eq:not02}
f(z) =\bigO\left(g(z)\right),  \quad z\in\calD,
\end{equation}
which means that there is  a constant $M$ such that $\vert f(z)\vert\le M \vert g(z)\vert$ for all $z\in\calD$.

Usually, for our asymptotic problems,  $z_0$ is the point of infinity, $\calD$ is an unbounded part of a sector, for example
\begin{equation}\label{eq:not03}
\calD=\left\{z\, \colon \vert z\vert \ge r , \  \alpha \le\phase\,z\le\beta\right\},
\end{equation}
where $\phase\,z$ denotes the phase of the complex number $z$, $r$ is a nonnegative number, and the real numbers $\alpha, \beta$ satisfy $\alpha\le\beta$.  When we write $f(z) =\bigO(1), z\in\calD$, we mean that $\vert f(z)\vert$ is bounded for all $z\in\calD$. 

The little $o$-symbol will also be used: 
\begin{equation}\label{eq:not04}
f(z) =o\left(g(z)\right),  \quad z\to z_0, \quad z\in\calD,
\end{equation}
which means that 
\begin{equation}\label{eq:not05}
\lim_{z\to z_0}
f(z)/g(z)=0,
\end{equation}
where the limiting point $z_0$ is approached inside $\calD$. When we write $f(z) =o(1), z\to z_0, z\in\calD$, we mean that $f(z)$ tends to zero when $z\to z_0$, $z\in\calD$. 

We write  
\begin{equation}\label{eq:not06}
f(z)\sim g(z),  \quad z\to z_0, \quad z\in\calD,
\end{equation}
when $\dsp{\lim_{z\to z_0} f(z)/g(z)=1}$, where the limiting point $z_0$ is approached inside $\calD$. In that case we say that the functions $f$ and $g$ are asymptotically equal at $z_0$.

\subsection{Asymptotic expansions}\label{sec:defexpas}    
Assume that, for $z\ne0$, a function $F$ has the representation 
\begin{equation}\label{ex:not07}
F(z)=a_0+\frac{a_1}z+\frac{a_2}{z^2}+\cdots+\frac{a_{n-
1}}{z^{n-1}}+R_n(z),\quad 
n=0,1,2,\ldots,
\end{equation}
with $F(z)=R_0(z)$, such that for each  $n =
0,1,2,\ldots$ the following relation holds
\begin{equation}\label{ex:not08}
R_n(z)=\bigO\left(z^{-n}\right),\quad {\rm as}\quad z\to\infty 
\end{equation}
in some
unbounded domain $\calD$. Then  $\dsp{\sn
a_nz^{-n}}$ is called an asymptotic expansion
of the
function $F(z)$, and we denote this by 
\begin{equation}\label{ex:not09}
F(z)\sim\sn a_nz^{-n},\quad z\to\infty, \quad z\in\calD.
\end{equation}
Poincar\'e introduced this form in 1886, and analogous forms 
can be given for
$z\to 0$, and so on.

Observe that we do not assume that the infinite series in \eqref{ex:not09}
converges for certain $z-$values. This
is not relevant in asymptotics; in the 
relation in \eqref{ex:not08} only a property of  $R_n(z)$ is requested, with
$n$ fixed.

\begin{example}[Exponential integral]\label{ex:exponint}
The classical example is the so-called 
exponential integral, that is,
\begin{equation}\label{eq:expint01}
F(z)=z\intp \frac{e^{-zt}}{t+1}\,dt=z\int_z^\infty
t^{-1}e^{z-t}\,dt=z\,e^zE_1(z),
\end{equation}
which we consider for $z>0$. Repeatedly using integration  by
parts in the second integral, we obtain
\begin{equation}\label{eq:expint02}
F(z)=1-\frac1z+\frac{2!}{z^2}-\cdots+\frac{(-1)^{n-1}(n-
1)!}{z^{n-1}}+
(-1)^nn!\,z\int_z^\infty\frac{e^{z-t}}{t^{n+1}}\,dt.
\end{equation}

In this
case we have, since $t\ge z$,
\begin{equation}\label{eq:expint03}
(-1)^nR_n(z)=n!\,z\int_z^\infty\frac{e^{z-t}}{t^{n+1}}\,dt\le
\frac{n!}{z^n}\int_z^\infty e^{z-t}\,dt=\frac{n!}{z^n}.
\end{equation}
Indeed, $R_n(z)=\bigO(z^{-n})$ as $z\to\infty$. Hence
\begin{equation}\label{eq:expint04}
z\int_z^\infty
t^{-1}e^{z-t}\,dt\sim\sn(-1)^n\frac{n!}{z^n},\quad z\to\infty.
\end{equation}
For extending this result to complex $z$, see Remark~\ref{rem:watlem}.
 \eoexample
\end{example}

\section{The classical methods for integrals}\label{sec:class}
We give a short overview of the classical methods. In each method the contributions to the asymptotic behavior of the integral can be obtained from one or a few {\em decisive points} of the interval of integration\footnote{Van der Corput used also the term {\em critical point}.}.

\subsection{Watson's lemma for Laplace-type integrals}\label{sec:aslap}
In this section we consider the large $z$ asymptotic expansions of Laplace-type integrals of the form
\begin{equation}\label{eq:watson01}
F_\lambda(z)=\intp t^{\lambda-1} f(t)e^{-zt}\,dt,\quad \Re z>0,  \quad \Re\lambda>0.
\end{equation}
For simplicity we assume that $f$  is analytic in a disc $\vert t\vert\le r$, $r>0$ and inside a sector $\calD\colon 
\alpha<\phase\,t<\beta$, where $\alpha<0$ and 
$\beta>0$. Also, we assume that there is a real number 
 $\sigma$ such that $f(t)=\bigO(e^{\sigma|t|})$ as
$t \to\infty$ in $\calD$. Then the substitution of $\dsp{f(t)= \sn a_n t^n}$ gives the asymptotic expansion 
\begin{equation}\label{eq:watson02}
F_\lambda(z)\sim\sn a_n\frac{\Gamma(n+\lambda)}{z^{n+\lambda}},\quad z\to\infty,
\end{equation}
which is valid inside the sector
\begin{equation}\label{eq:watson03}
-\beta-\tfrac12\pi+\delta\le\phase\,z\le-\alpha+\tfrac12\pi-\delta,
\end{equation}
where $\delta$ is a small positive number.

This useful result is known as Watson's lemma. By using a finite expansion of $f$ with remainder we can prove the asymptotic property needed for the Poincar\'e-type expansion. For the proof (also for more general conditions on $f$) we refer to \cite[p.~114]{Olver:1997:ASF}.

To explain how the bounds in \eqref{eq:watson03} arise, we allow the path of integration
in \eqref{eq:watson01} to turn over an angle $\tau$, 
and write $\phase\,t=\tau$, 
$\phase\,z=\theta$, where $\alpha<\tau<\beta$. The 
condition for convergence in
\eqref{eq:watson01} is $\cos(\tau+\theta)>0$, that is, $-
\frac12\pi<\tau+\theta<\frac12\pi$.
Combining this with  the bounds for $\tau$ we obtain the bounds 
for
$\theta$ in \eqref{eq:watson03}.

Several modifications of Watson's lemma are considered in the literature, for example of Laplace transforms with logarithmic singularities at the origin.

\begin{remark}\label{rem:watlem}
In Example~\ref{ex:exponint} for the exponential integral we have  
$f(t)=(1+t)^{-1}$. In that case $f$ is analytic in the 
sector $|\phase\,t|<\pi$, and  
$\alpha=-\pi,\,
\beta=\pi$. Hence, the asymptotic expansion given in \eqref{eq:expint04}
holds in the sector $|\phase\,z|\le\frac32\pi-\delta$. This range is
much larger than the usual domain of definition for the
exponential integral, which reads: $|\phase\,z|<\pi$. 
\eoremark
\end{remark}

\begin{example}\label{ex:modbes}
The modified Bessel function $K_\nu(z)$ has the integral representation
\begin{equation}\label{eq:modb01}
\mathop{K_{{\nu}}\/}\nolimits\!\left(z\right)=
\frac{\pi^{{\frac{1}{2}}}(2z)^{\nu}e^{-z}}{\mathop{\Gamma\/}\nolimits\!\left(\nu+\frac{1}{2}\right)}\int _{0}^{\infty}t^{\nu-\frac12}e^{{-2zt}}f(t) \,dt,
\quad \Re\nu>-\tfrac12,
\end{equation}
where
\begin{equation}\label{eq:modb02}
f(t)=(t+1)^{{\nu-\frac{1}{2}}}=\sn \binom{\nu-\frac12}{n}t^n,\quad \vert t\vert <1.
\end{equation}
Watson's lemma can be used by substituting this expansion into \eqref{eq:modb01}, and we obtain
\begin{equation}\label{eq:modb03}
K_\nu(z)\sim\sqrt{{\frac\pi{2z}}}e^{-z}\sn\frac{c_n(\nu)}{n!\,(8z)^{n}},\quad |\phase\, z
|\le\tfrac32\pi-\delta,
\end{equation}
where $\delta$ is a small positive number. The coefficients are given by $c_0(\nu)=1$ and 
\begin{equation}\label{eq:modb04}
c_n(\nu)=\left(4\nu^2-1\right)\left(4\nu^2-3^2\right)\cdots\left(4\nu^2-(2n-1)^2\right),\quad n\ge1.
\end{equation}

The expansion in \eqref{eq:modb03} is valid for bounded values of $\nu$. Large values of $\nu$ destroy the asymptotic nature and the expansion in \eqref{eq:modb03} is valid for large $z$ in the shown sector, uniformly with respect to bounded values of $\nu$. 

The expansion terminates when $\nu=\pm\frac12, \pm\frac32,\ldots$, and we have a finite exact representation. By using the many relations between the Bessel functions, expansions for large argument of all other Bessel functions can be derived  from this example. For all these expansions sharp estimates are available of remainders in the expansions; see Olver's work \cite[Chapter~7]{Olver:1997:ASF}, which is based on using differential equations.

There are many integral representations for Bessel functions. For $K_\nu(z)$ we have, for example,
\begin{equation}\label{eq:modb05}
K_\nu(z)=\intp e^{-z\cosh t}\cosh(\nu t)\,dt,\quad \Re z>0.
\end{equation}
It is not difficult to transform this integral into the standard form given in \eqref{eq:watson01}, for example by substituting $\cosh t=1+s$, or to consider it as an example for Laplace's method considered in \S\ref{sec:laplacemethod}. However, then the method for obtaining the explicit form of the coefficients as shown in \eqref{eq:modb04} will be rather complicated.
\eoexample
\end{example}

\subsection{The method of Laplace}\label{sec:laplacemethod}

In Laplace's method we deal with integrals of the form
\begin{equation}\label{eq:Laplace00}
F(z)=\int_a^b e^{-z p(t)} q(t)\,dt,
\end{equation}
(see \cite[\S3.7]{Olver:1997:ASF}), where it is assumed that these integrals can be transformed into
integrals of the form
\begin{equation}\label{eq:Laplace01}
F(z)=\intr e^{-z t^2} f(t)\,dt,\quad \Re z>0.
\end{equation}
The function $f$ is assumed to be analytic inside a domain $\calD$ of the complex plane that contains the real axis in its interior.
 
By splitting the contour of integration in two parts, for positive and negative $t$, two integrals arise that can be expanded by applying Watson's lemma considered in \S\ref{sec:aslap}.
 
On the other hand we can substitute a Maclaurin expansion $\dsp{f(t)=\sk c_k t^k}$    to obtain
\begin{equation}\label{eq:Laplace02}
F(z)\sim \sqrt{\frac{\pi}{z}}\sk \left(\tfrac12\right)_k\frac{c_{2k}}{z^k}, \quad z\to\infty.
\end{equation}
The domain of validity depends on the location of the singularities of the function $f$ in the complex plane. We can suppose that this function is even.

As in Watson's lemma we assume that $\calD$  contains a disk $\vert t\vert\le r$, $r>0$, and  sectors $\alpha<\phase(\pm t)<\beta$, where $\alpha<0$ and $\beta>0$.
In addition we need convergence conditions, say by requiring that there is  a real number $\sigma$ such that $f(t)=\bigO(e^{\sigma|t^2|})$ as $t \to\infty$ in $\calD$. 

Write $z=r e^{i\theta}$ and $t=\sigma e^{i\tau}$. Then, when rotating the path, the condition for convergence at infinity is $\cos(\theta+2\tau)>0$, that is,  $-\frac12\pi<\theta+2\tau<\frac12\pi$, which should be combined with $\alpha<\tau<\beta$ for staying inside the sector. Then the expansion in \eqref{eq:Laplace02} holds uniformly  inside the sector
\begin{equation}\label{eq:Laplace03}
-2\beta-\tfrac12\pi+\delta\le\phase\,z\le\tfrac12\pi-2\alpha-\delta,
\end{equation}
for any small positive number $\delta$.

 For a detailed analysis of Laplace's method we refer to \cite[pp.~121--127]{Olver:1997:ASF}.

\subsection{The saddle point method and paths of steepest descent}\label{sec:sad}

In this case the integrals are presented as contour integrals in the complex plane:
\begin{equation}\label{eq:saddle01}
F(z)=\int_\calC e^{-z\phi(w)}\psi(w)\,dw,
\end{equation}
where $z$ is a large real or complex parameter. The functions  $\phi, \psi$ are  analytic in a domain $\calD$ of the 
complex plane. The integral is taken along a path $\calC$ in $\calD$,  
and avoids the singularities and branches of the integrand. 
Integrals of this type arise naturally in the context  of linear wave 
propagation and in other physical problems;  many special functions can be represented by such integrals.

Saddle points are zeros of $\phi^\prime(w)$ and for the integral we select modifications of the contour $\calC$, usually through a saddle point and on the path we require that $\Im(z\phi(w))$ is constant, a steepest descent path. After several steps we may obtain representations as in \eqref{eq:Laplace01}, after which we can apply Laplace's method. For an extensive discussion of the saddle point method we refer to De Bruijn's book \cite{Bruijn:1981:AMA}.

\subsection{Generating functions; Darboux's method}\label{sec:darb}
The classical orthogonal polynomials, and many other special functions, have  generating functions of the form\begin{equation}\label{eq:gfin01}
G(z,w)=\sn F_n(z) w^n.
\end{equation}
The radius of convergence may be finite or infinite, and may depend on the variable $z$. For example, the Laguerre polynomials satisfy the relation
\begin{equation}\label{eq:gfin02}
(1-w)^{-\alpha-1}e^{-wz/(1-w)}=\sn L_n^\alpha(z)w^n,\quad \vert w\vert <1;
\end{equation}
$\alpha$ and $z$ may assume any finite complex value.

From the generating function a representation in the form a Cauchy integral follows:
\begin{equation}\label{eq:gfin03}
F_n(z)=\frac{1}{2\pi i}\int_{\calC} G(z,w)\,\frac{dw}{w^{n+1}},
\end{equation}
where $\calC$ is a circle around the origin inside the domain where $G(z,w)$ is  analytic as a function of $w$. 

When the function $G(z,w)$ has simple algebraic singularities, an asymptotic expansion of $F_n(z)$ can usually  be obtained by deforming the contour $\calC$ around the branch points or other singularities of  $G(z,w)$ in the $w-$plane. 

\begin{example}[Legendre polynomials]\label{ex:legenpol}
The Legendre polynomials have the generating function
\begin{equation}\label{eq:gfin04}
\frac{1}{\sqrt{1-2xw+w^2}}=\sn P_n(x)w^n,\quad -1\le x\le
1,\quad |w|<1.
\end{equation}
We write $x=\cos\theta$, $0\le \theta\le \pi$,  and obtain the Cauchy integral representation
\begin{equation}\label{eq:gfin05}
P_n(\cos\theta)=\frac{1}{2\pi i}\int_\calC \frac{1}{\sqrt{1-2\cos\theta \,w+w^2}}\frac{dw}{w^{n+1}},
\end{equation}
where $\calC$ is a circle around the origin with radius less than $1$. There are two singular points on the unit circle: $w_\pm= e^{\pm i\theta}$.  When $\theta\in(0,\pi)$ we can deform the contour $\calC$ into two loops $\calC_\pm$ around the two branch cuts. Because $P_n(-x)=(-1)^nP_n(x)$, we take $0\le x<1$, $0<\theta\le\frac12\pi$. 

We assume that the square root in \eqref{eq:gfin05}  is positive for real values of $w$ and that  branch cuts run from each $w_\pm$ parallel to the real axis, with $\Re w\to +\infty$. For $\calC_+$ we substitute $w=w_+e^s$, and obtain a similar contour  $\calC_+$ around the origin in the $s-$plane. We start the integration along $\calC_+$ at $+\infty$, with $\phase\,s=0$, turn around the origin in the clock-wise direction, and return to $+\infty$ with $\phase\,s=-2\pi$. The contribution from $\calC_+$ becomes
\begin{equation}\label{eq:gfin06}
P_n^+(\cos\theta)=\frac{e^{-\left(n+\frac12\right)i\theta+\frac14\pi i}}{\pi\sqrt{2\sin\theta}}\intp e^{-ns}f^+(s)\,\frac{ds}{\sqrt{s}},
\end{equation}
where
\begin{equation}\label{eq:gfin07}
f^+(s)=\sqrt{\frac{s}{e^s-1}\,\frac{1-e^{-2i\theta}}{e^{s}-e^{-2i\theta}}}.
\end{equation}
Expanding $f^+$ in powers of $s$, we can use Watson's lemma (see \S\ref{sec:aslap}) to obtain the large $n$ expansion of $P_n^+(\cos\theta)$. If $\theta\in[\theta_0,\frac12\pi]$, where $\theta_0$ is a small positive number, then the conditions are satisfied to apply Watson's lemma. 

The contribution from the singularity at $w_-$ can be obtained in the same way. It is the complex conjugate of the contribution from $w_+$, and we have $P_n(\cos\theta)=2\Re P_n^+(\cos\theta)$.
\eoexample
\end{example}

For an application to large degree asymptotic expansions of generalized  Bernoulli and Euler polynomials, see \cite{Lopez:2010:LDA}; see also  \cite[Chapter~2]{Wong:2001:AAI}. 

The way of handling coefficients of power series is related to  {\em Darboux's method}, in which again the asymptotic behavior is considered of the coefficients of a power series $f(z)=\sum a_n z^n$. A comparison function $g$, say, is needed  with the same relevant singular point(s) as $f$. When $g$ has an expansion $g(z)=\sum b_n z^n$, in which the coefficients $b_n$ have known asymptotic behavior, then, under certain conditions on $f(z)-g(z)$ near the singularity, it is possible to find asymptotic forms for the coefficients $a_n$. For an introduction to Darboux's method and examples for orthogonal polynomials, we refer to \cite[\S8.4]{Szego:1975:OP}.

For the Laguerre polynomials the method as explained in Example~\ref{ex:legenpol} does not work. The essential singularity at $w=1$ in the left-hand side of \eqref{eq:gfin02} requires a different approach. See \cite{Perron:1921:UDV}, where the more general case of Kummer functions is considered and, more recently, \cite{Borwein:2008:ELA} for Laguerre polynomials $L_n^\alpha(z)$ for large values of $n$, with complex $z\notin[0,\infty)$.

When relevant singularities are in close proximity, or even coalescing, we need uniform methods and for  a uniform treatment  we refer to \cite{Wong:2005:OAU}. In our example of the Legendre polynomials, uniform methods are needed  to deal with small values of $\theta$; in that case $J-$Bessel functions are needed.

\subsection{Mellin-type integrals}\label{sec:mellin}
Mellin convolution integrals are of the form 
\begin{equation}\label{eq:Mell01}
F_{\lambda}(x)=\intp t^{\lambda-1} h(xt) f(t)\,dt,
\end{equation}
and they reduce to the standard form of Watson's lemma when $h(t)=e^{-t}$. For the general case \eqref{eq:Mell01} powerful methods have been developed, in which asymptotic expansions are obtained for $x\to0$ and for $x\to\infty$. For more details we refer to \cite{Lopez:2008:AMC}, \cite[Chapters 4,6]{Bleistein:1975:AEI} and  \cite[Chapter~VI]{Wong:2001:AAI}.

A main step in the method to obtain asymptotic expansions of the integral in \eqref{eq:Mell01} is the use of Mellin transforms and their inverses. These inverses can be viewed as Mellin-Barnes integrals. An example is the Meijer $G-$function
\begin{equation}\label{eq:Mell02}
\begin{array}{ll}
\dsp{\mathop{{G^{{m,n}}_{{p,q}}}\/}\nolimits\!\left(z;{a_{1},\dots,a_{p}\atop b_{1},\dots,b_{q}}\right)}=\\[8pt]
\quad\quad
\dsp{\frac{1}{2\pi i}\mathlarger{\int}_{{\calL}}{\frac{\prod\limits _{{\ell=1}}^{m}\mathop{\Gamma\/}\nolimits\!\left(b_{{\ell}}-s\right)\prod\limits _{{\ell=1}}^{n}\mathop{\Gamma\/}\nolimits\!\left(1-a_{{\ell}}+s\right)}{\prod\limits _{{\ell=m}}^{{q-1}}\mathop{\Gamma\/}\nolimits\!\left(1-b_{{\ell+1}}+s\right)\prod\limits _{{\ell=n}}^{{p-1}}\mathop{\Gamma\/}\nolimits\!\left(a_{{\ell+1}}-s\right)}} z^{s}\, ds.}
\end{array}
\end{equation}
It can be viewed as the inversion of a Mellin transform. Many special functions can be written in terms of this function.  
The integration path $\calL$ separates the poles of the factors  $\Gamma(b_\ell-s)$ from those of the factors $\Gamma(1-a_\ell+s)$. 
By shifting the contour and picking up the residues, expansions for $z\to0$ (when shifting to the right) and $z\to\infty$  (when shifting to the left) may be obtained. For an extensive treatment we refer to \cite{Paris:2001:AMB}.

\section{The method of stationary phase}\label{sec:stat}
We consider this method in more detail, and we give also new elements which will give uniform expansions.
The integrals are of the type
\begin{equation}\label{eq:mstph01}
F(\omega)=\int_a^b e^{i\omega\phi(t)}\psi(t)\,dt,
\end{equation}
where $\omega$ is a real large parameter, $a,b$ and $\phi$ are real; $a=-\infty$ or/and $b=+\infty$ are allowed.
The idea of the method of stationary phase is originally developed by Stokes and Kelvin. 

The asymptotic character of the integral in \eqref{eq:mstph01} is completely determined if the behavior of the functions $\phi$ and $\psi$ is known in the vicinity of the decisive points. These are
\begin{itemize}
\item
stationary points: zeros of $\phi^{\prime}$ in $[a,b]$;
\item
the finite endpoints $a$ and $b$;
\item
values in or near $[a,b]$ for which $\phi$ or $\psi$ are singular.
\end{itemize}
The following integral shows all these types of decisive points:
\begin{equation}\label{eq:mstph02}
F(\omega)=\int_{-1}^1 e^{i\omega t^2}\sqrt{\left\vert t-c\right\vert}\,dt,\quad -1<c<1.
\end{equation}

\begin{figure}
\begin{center}
\epsfxsize=6cm \epsfbox{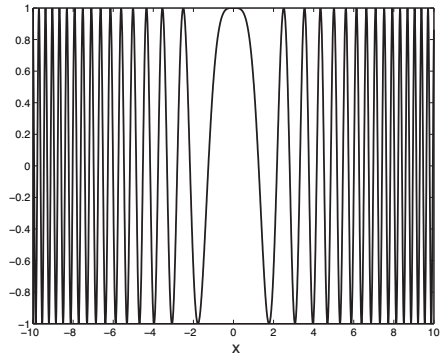}
\end{center}
 \caption{{\small $\Re\left(e^{i\omega x^2}\right)$ for $\omega=1$ with at the stationary point $x=0$ less oscillations. \label{fig:StPhase}}}
\end{figure}

In Figure~\ref{fig:StPhase} we see why a stationary point may give a contribution: less oscillations occur at the stationary point compared with other points in the interval, where the oscillations neutralize each other.

\subsection{Integrating by parts: no stationary points}\label{sec:stphinpa}

If in \eqref{eq:mstph01}  $\phi$ has no stationary point in the finite interval $[a,b]$,  and $\phi$, $\psi$ are regular on $[a,b]$, then contributions from the endpoints follow from integrating by parts.  We have
\begin{equation}\label{eq:mstph03}
\begin{array}{@{}r@{\;}c@{\;}l@{}}
F(\omega)&=&
\dsp{\int_a^b e^{i\omega\phi(t)}\psi(t)\,dt
=\frac{1}{i\omega} \int_{a}^b \psi(t)\frac{de^{i\omega\phi(t)}}{\phi^\prime(t)} =}\\[8pt]
&=&
\dsp{\frac{e^{i\omega\phi(b)}}{i\omega\phi^\prime(b)}\psi(b)-\frac{e^{i\omega\phi(a)}}{i\omega\phi^\prime(a)}\psi(a)+\frac{1}{i\omega}\int_{a}^b e^{i\omega\phi(t)}\psi_1(t)\,dt,}
\end{array}
\end{equation}
where
\begin{equation}\label{eq:mstph04}
\psi_1(t)=-\frac{d}{dt}\,\frac{\psi(t)}{\phi^\prime(t)}.
\end{equation}
The integral in this result has the same form as the one in \eqref{eq:mstph01}, and when $\phi$ and $\psi$ are sufficiently smooth, we can continue this procedure.

In this way we obtain for $N=0,1,2,\ldots$ the compound expansion 
\begin{equation}\label{eq:mstph05}
\begin{array}{@{}r@{\;}c@{\;}l@{}}
F(\omega)&=&
\dsp{\frac{e^{i\omega\phi(b)}}{\phi^\prime(b)}\sum_{n=0}^{N-1}\frac{\psi_n(b)}{(i\omega)^{n+1}}-\frac{e^{i\omega\phi(a)}}{\phi^\prime(a)}\sum_{n=0}^{N-1}\frac{\psi_n(a)}{(i\omega)^{n+1}}}\ +\\[8pt]
&&\dsp{\frac{1}{(i\omega)^N}\int_{a}^b e^{i\omega\phi(t)}\psi_N(t)\,dt,}
\end{array}
\end{equation}
where for $N=0$ the sums in first line vanish. The integral can be viewed as a remainder of the expansion. We have $\psi_0=\psi$ and 
\begin{equation}\label{eq:mstph06}
\psi_{n+1}(t)=-\frac{d}{dt}\,\frac{\psi_n(t)}{\phi^\prime(t)},\quad n=0,1,2,\ldots.
\end{equation}

This expansion can be obtained when $\phi,\psi\in C^N[a,b]$. When we assume that we can find positive numbers $M_N$ such that
$\left\vert \psi_N(t)\right\vert\le M_N$ for $t\in[a,b]$, we can find an upper bound of the remainder in \eqref{eq:mstph05}, and this estimate will be of order $\bigO\left(\omega^{-N}\right)$. Observe that the final term in each of the sums in \eqref{eq:mstph05} are
also of order $\bigO\left(\omega^{-N}\right)$.

\subsection{Stationary points and the use of neutralizers}\label{sec:stphneut}
When $\phi$ in \eqref{eq:mstph01} has exactly one simple internal stationary point, say at $t=0$, that is, $\phi^\prime(0)=0$ and $\phi^{\prime\prime}(0)\ne0$, we can transform the integral into the form
\begin{equation}\label{eq:mstph07}
F(\omega)=\int_{a}^b e^{i\omega t^2}f(t)\,dt,\quad \omega>0,
\end{equation}
possibly with different $a$ and $b$ as in \eqref{eq:mstph01}, and again we assume finite $a<0$ and $b>0$. Because of the stationary point at the origin, the straightforward integration by parts method of the previous section cannot be used. 

We assume that $f\in C^N[a,b]$ for some positive $N$. There are three decisive points: $a, 0, b$, and we can split up the interval into $[a,0]$ and $[0,b]$, each subinterval having two decisive points. To handle this, Van der Corput \cite{Corput:1948:MCP} introduced neutralizers in order to get  intervals in which only one decisive point exists.

A  {\em neutralizer} $N_c$ at a point $c$  is a $C^\infty(\RR)$ function such that:
\begin{enumerate}
\item
$N_c(c)=1$, and all its derivatives vanish at $c$.

\item
There is a positive number $d$ such that $N_c(x)=0$ outside $(c-d,c+d)$.
\end{enumerate}
It is not difficult to give explicit forms of such a neutralizer, but in the analysis this is not needed.

With a neutralizer at the origin we can write the integral in \eqref{eq:mstph07} in the form
\begin{equation}\label{eq:mstph08}
F(\omega)=
\intr e^{i\omega t^2}f(t) N_0(t)\,dt+
\int_{a}^b e^{i\omega t^2}f(t)(1-N_0(t))\,dt.
\end{equation}

We write in the first integral 
\begin{equation}\label{eq:mstph09}
f(t)N_0(t)=\sum_{n=0}^{N-1} c_nt^n +R_N(t)
\end{equation}
where the coefficients $c_n=f^{(n)}(0)/n!$ do not depend on the neutralizer, because all derivatives of $N_0(t)$ at $t=0$ vanish. 

We would like to evaluate the integrals by extending the interval $[a,b]$ to $\RR$. That is, we evaluate
\begin{equation}\label{eq:mstph10}
\intr e^{i\omega t^2} t^{n}\,dt
\end{equation}
of which the ones with odd $n$ vanish. The other ones with $n\ge2$ diverge. 

For a way out we refer to \S\ref{sec:stphavneut}. Other methods are based on introducing a converging factor, such as $e^{-\eps t^2}$, $\eps>0$, in the integral in \eqref{eq:mstph07}. For details on this method, see \cite{Olver:1974:EBS}. In \cite[\S II.3]{Wong:2001:AAI} a method of Erd\'elyi \cite{Erdelyi:1955:ARF} is used.

In the second integral in \eqref{eq:mstph08}, the stationary point at $t=0$ is harmless (it is neutralized), and we can integrate by parts to obtain the contributions from $a$ and $b$, as we have done in \S\ref{sec:stphinpa}. The expansion of  the second integral is the same as in \eqref{eq:mstph05}, now  with $\phi(t)=t^2$. Again, the asymptotic terms do not depend on the neutralizer $N_0$. 

\subsection{How to avoid neutralizers?}\label{sec:stphavneut}
Neutralizers can be used in more complicated problems, and Van der Corput's work on the method of stationary phase was pioneering. After 1960 the interest in uniform expansions came up and he once admitted (private communications, 1967) that he didn't see how to use his neutralizers when two or more decisive points are nearby or even coalescing. This became a major topic in uniform methods, as we will see in later sections.

Bleistein \cite{Bleistein:1966:UAE} introduced in 1966 a new method on integrating by parts in which two decisive points were allowed to coalesce.\footnote{Van der Corput wrote a review of this paper in Mathematical Reviews.} The method was not used for the method of stationary phase, but as we will see we can use a simple form of Bleistein's method also for this case. 

Again we consider  \eqref{eq:mstph07} with finite $a<0$ and $b>0$ and $f\in C^N[a,b]$.
We write $f(t)=f(0)+(f(t)-f(0))$ and obtain
\begin{equation}\label{eq:mstph12}
F(\omega)=f(0)\Phi_{a,b}(\omega)+\int_{a}^b e^{i\omega t^2}\left(f(t)-f(0)\right)\,dt,
\end{equation}
where
\begin{equation}\label{eq:mstph13}
\Phi_{a,b}(\omega)=\int_{a}^b e^{i\omega t^2}\,dt,
\end{equation}
a Fresnel-type integral \cite{Temme:2010:ERF}. 

Because the integrand vanishes at the origin we can integrate by parts in \eqref{eq:mstph12}, and write 
\begin{equation}\label{eq:mstph14}
F(\omega)=f(0)\Phi_{a,b}(\omega)+\frac{1}{2i\omega}\int_{a}^b \frac{f(t)-f(0)}{t}\,de^{i\omega t^2}.
\end{equation}
This gives
\begin{equation}\label{eq:mstph15}
\begin{array}{@{}r@{\;}c@{\;}l@{}}
F(\omega)&=&\dsp{\frac{e^{i\omega b^2}}{2i\omega}\frac{f(b)-f(0)}{b}-
\frac{e^{i\omega a^2}}{2i\omega}\frac{f(a)-f(0)}{a}}\ +\\[8pt]
&&\dsp{f(0)\Phi_{a,b}(\omega)+\frac{1}{2i\omega}\int_{a}^b f_1(t) e^{i\omega t^2}\,dt,}
\end{array}
\end{equation}
where 
\begin{equation}\label{eq:mstph16}
f_1(t)=-\frac{d}{dt}\frac{f(t)-f(0)}{t}.
\end{equation}

The integral in \eqref{eq:mstph15} can be expanded in the same way. We obtain
\begin{equation}\label{eq:mstph17}
\begin{array}{@{}r@{\;}c@{\;}l@{}}
F(\omega)&=&\dsp{
{e^{i\omega b^2}}{}\sum_{n=0}^{N-1}\frac{C_n(b)}{(2i\omega)^{n+1}}-
{e^{i\omega a^2}}{}\sum_{n=0}^{N-1}\frac{C_n(a)}{(2i\omega)^{n+1}}}\ +\\[8pt]
&&\dsp{\Phi_{a,b}(\omega)\sum_{n=0}^{N-1}\frac{f_n(0)}{(2i\omega)^n}+\frac{1}{(2i\omega)^N}\int_{a}^b f_N(t) e^{i\omega t^2}\,dt,}
\end{array}
\end{equation}
where $N\ge0$ and for $n=0,1,2,\ldots,N$ we define
\begin{equation}\label{eq:mstph18}
C_n(t)=\frac{f_n(t)-f_n(0)}{t}, \quad
f_{n+1}(t)=-\frac{d}{dt}\frac{f_n(t)-f_n(0)}{t}.
\end{equation}
For $N=0$ the terms with the sums in \eqref{eq:mstph17} all vanish.

The integral in  \eqref{eq:mstph17} can be viewed as a remainder of the expansion. When we assume that we can find positive numbers $M_N$ such that
$\left\vert f_N(t)\right\vert\le M_N$ for $t\in[a,b]$, we can find an upper bound of this remainder, which will be of order $\bigO\left(\omega^{-N}\right)$. The final term in each of the sums in \eqref{eq:mstph17} are
also of order $\bigO\left(\omega^{-N}\right)$.

If we wish we can keep the Fresnel-type integral in the expansion in \eqref{eq:mstph17} as it is, but we can also proceed by writing
\begin{equation}\label{eq:mstph19}
\begin{array}{@{}r@{\;}c@{\;}l@{}}
\Phi_{a,b}(\omega)
&=&\dsp{\intr e^{i\omega t^2}\,dt - \int_{-\infty}^a e^{i\omega t^2}\,dt - \int_{b}^\infty e^{i\omega t^2}\,dt}\\[8pt]
&=&\dsp{\sqrt{\frac{\pi}{\omega}}e^{\frac14\pi i}- \int_{-a}^\infty e^{i\omega t^2}\,dt - \int_{b}^\infty e^{i\omega t^2}\,dt}.
\end{array}
\end{equation}
The integrals in the second line can be expressed in terms of the complementary error function (see \eqref{eq:waerden02}) and can be expanded for large values of $\omega$ by using \eqref{eq:bruijn09} (when $a$ and $b$ are bounded away from zero); see \cite[\S7.12(ii)]{Temme:2010:ERF}. In this way we can obtain, after some re-arrangements, the same expansions as with the neutralizer method of \S\ref{sec:stphneut}. However,
we observe the following points.

\begin{itemize}
\item
When we do not expand the function $\Phi_{a,b}(\omega)$ for large $\omega$, the expansion in  \eqref{eq:mstph17} remains valid when $a$ and/or $b$ tend to zero. This means, we can allow the decisive points to coalesce (Van der Corput's frustration), and  in fact we have an expansion that is valid uniformly with respect to small values of $a$ and $b$.
\item
There is no need to handle the divergent integrals in \eqref{eq:mstph10}.  As we explained after  \eqref{eq:mstph10}, a rigorous approach is possible by using converging factors, or other methods, but the present representation is very convenient.

\item
The remainders in the neutralizer method depend on the chosen neutralizers; the remainder in \eqref{eq:mstph17} only depends on $f$ and its derivatives.
\end{itemize}

\subsection{Algebraic singularities at both endpoints}\label{sec:singend}
We assume that $f$ is $N$ times continuously differentiable in the finite interval $[\alpha,\beta]$. Consider the following integral
\begin{equation}\label{eq:mstph20}
F(\omega)=\int_{\alpha}^{\beta}e^{i\omega t}(t-\alpha)^{\lambda-1}(\beta-t)^{\mu-1} f(t)\,dt,
\end{equation}
where $\Re \lambda>0$, $\Re \mu >0$. A straightforward approach using integrating by parts is not possible. However, see \S\ref{sec:singendpart}.

For this class of integrals we have the following result 
\begin{equation}\label{eq:mstph21}
F(\omega)=A_N(\omega)+B_N(\omega)+\bigO\left(1/\omega^N\right),\quad \omega\to\infty,
\end{equation}
where
\begin{equation}\label{eq:mstph22}
\begin{array}{ll}
\dsp{A_N(\omega)=\sum_{n=0}^{N-1}\frac{\Gamma(n+\lambda)}{n!\,\omega^{n+\lambda}}e^{i(\frac12\pi(n+\lambda)+\alpha\omega)}\left.\frac{d^n}{dt^n}\left((\beta-t)^{\mu-1} f(t)\right)\right\vert_{t=\alpha},}\\[8pt]
\dsp{B_N(\omega)=\sum_{n=0}^{N-1}\frac{\Gamma(n+\mu)}{n!\,\omega^{n+\mu}}e^{i(\frac12\pi(n-\mu)+\beta\omega)}\left.\frac{d^n}{dt^n}\left((t-\alpha)^{\lambda-1} f(t)\right)\right\vert_{t=\beta}.}
\end{array}
\end{equation}
Erd\'elyi's proof in \cite{Erdelyi:1955:ARF} is based on the use of neutralizers and the result can be applied to, for example,
the Kummer function written in the form
\begin{equation}\label{eq:mstph24}
\CHF{\lambda}{\lambda+\mu}{i\omega}=\frac{\Gamma(\lambda+\mu)}{\Gamma(\lambda)\Gamma(\mu)}
\int_{0}^{1}e^{i\omega t}t^{\lambda-1}(1-t)^{\mu-1} \,dt, \quad \Re\lambda, \Re \mu >0.
\end{equation}
By using \eqref{eq:mstph21} a standard result for this function with $\omega\to+\infty$ follows, see \cite[\S13.7(i)]{Olde:2010:CHF}, but there are more elegant ways to obtain a result valid for general complex argument.

\subsubsection{Integrating by parts}\label{sec:singendpart}
Also in this case we can avoid the use of neutralizers by using integration by parts. This cannot be done in a straightforward way because of the singularities of the integrand in \eqref{eq:mstph20} at the endpoints. Observe that Erd\'elyi's expansion gives two expansions, one from each decisive endpoint.  By a modification of Bleistein's procedure we obtain expansions which take contributions from both decisive points in each step.  

We write (see also the method used in \S\ref{sec:stphavneut})
\begin{equation}\label{eq:mstph25}
f(t)=a_0+b_0(t-\alpha)+(t-\alpha)(\beta-t)g_0(t),
\end{equation}
where $a_0, b_0$ follow from substituting $t=\alpha$ and $t=\beta$. This gives 
\begin{equation}\label{eq:mstph26}
a_0=f(\alpha),\quad b_0=\frac{f(\beta)-f(\alpha)}{\beta-\alpha},
\end{equation}
and for \eqref{eq:mstph20} we obtain
\begin{equation}\label{eq:mstph27}
F(\omega)=a_0\Phi+b_0\Psi+\int_{\alpha}^{\beta}e^{i\omega t}(t-\alpha)^{\lambda}(\beta-t)^{\mu} g_0(t)\,dt,
\end{equation}
where (see \eqref{eq:mstph24})
\begin{equation}\label{eq:mstph28}
\begin{array}{ll}
\Phi=(\beta-\alpha)^{\lambda+\mu-1}e^{i\omega\alpha}\dsp{\frac{\Gamma(\lambda)\Gamma(\mu)}{\Gamma(\lambda+\mu)}}\CHF{\lambda}{\lambda+\mu}{i(\beta-\alpha)\omega},\\[8pt]
\Psi=(\beta-\alpha)^{\lambda+\mu}e^{i\omega\alpha}\dsp{\frac{\Gamma(\lambda+1)\Gamma(\mu)}{\Gamma(\lambda+\mu+1)}}\CHF{\lambda+1}{\lambda+\mu+1}{i(\beta-\alpha)\omega}.
\end{array}
\end{equation}
Now we integrate by parts and obtain, observing that the integrated terms will vanish,
\begin{equation}\label{eq:mstph29}
F(\omega)=a_0\Phi+b_0\Psi+\frac{1}{i\omega}\int_{\alpha}^{\beta}e^{i\omega t}(t-\alpha)^{\lambda-1}(\beta-t)^{\mu-1} f_1(t)\,dt,
\end{equation}
where
\begin{equation}\label{eq:mstph30}
f_1(t)=-(t-\alpha)^{1-\lambda}(\beta-t)^{1-\mu}\frac{d}{dt}\left((t-\alpha)^{\lambda}(\beta-t)^{\mu} g_0(t)\right).
\end{equation}

We can continue with this integral in the same manner, and obtain
\begin{equation}\label{eq:mstph31}
F(\omega)=
\Phi\sum_{n=0}^{N-1}\frac{a_n}{(i\omega)^n}+
\Psi\sum_{n=0}^{N-1}\frac{b_n}{(i\omega)^n}+R_N(\omega),\quad N=0,1,2,\ldots,
\end{equation}
where
\begin{equation}\label{eq:mstph32}
R_N(\omega)=
\frac{1}{(i\omega)^N}\int_{\alpha}^{\beta}e^{i\omega t}(t-\alpha)^{\lambda-1}(\beta-t)^{\mu-1} f_N(t)\,dt.
\end{equation}
The coefficients are defined by
\begin{equation}\label{eq:mstph33}
a_n=f_n(\alpha),\quad b_n=\frac{f_n(\beta)-f_n(\alpha)}{\beta-\alpha},
\end{equation}
and the functions $f_n$ follow from the recursive scheme
\begin{equation}\label{eq:mstph34}
\begin{array}{@{}r@{\;}c@{\;}l@{}}
f_n(t)&=&a_n+b_n(t-\alpha)+(t-\alpha)(\beta-t)g_n(t),\\[8pt]
f_{n+1}(t)&=&\dsp{-(t-\alpha)^{1-\lambda}(\beta-t)^{1-\mu}\frac{d}{dt}\left((t-\alpha)^{\lambda}(\beta-t)^{\mu} g_n(t)\right)},
\end{array}
\end{equation}
with $f_0=f$.

The expansion in \eqref{eq:mstph31} contains confluent hypergeometric functions, 
and Erd\'elyi's 
expansion in \eqref{eq:mstph21} is in terms of elementary functions. 
Erd\'elyi's expansion follows by using the asymptotic expansion of the confluent hypergeometric function. Estimates of $R_N(\omega)$ follow when we know more about the function $f$, in particular when we have bounds on the derivatives of $f$. 

Again, the present approach avoids neutralizers, and the expansion remains valid when the decisive points $\alpha$ and $\beta$ coalesce. Of course, when $\omega(\beta-\alpha)$ is not large,  the ${}_1F_1-$functions in \eqref{eq:mstph28} should not be expanded.

\subsection{Some examples where the method fails}\label{sec:stphaseex}

In some cases it is not possible to obtain quantitative
information in the method of stationary phase. For instance, when $\phi\in C^\infty[a,b]$ does not have a stationary 
point
in $[a,b]$, and $\psi\in C^\infty[a,b]$ vanishes with all its derivatives at $a$ and
$b$. In that case we have
\begin{equation}\label{eq:mstph35}
F(\omega)=\bigO\left(\omega^{-n}\right)
\end{equation}
as $\omega\to\infty$, for any $n$, and we say the function $F(\omega)$ is exponentially small. An example is
\begin{equation}\label{eq:mstph36}
F(\omega)=\intr e^{i\omega t}\frac{dt}{1+t^2}= \pi e^{-\omega}.
\end{equation}
The integral
\begin{equation}\label{eq:mstph37}
G(\omega)=\int_{-\frac12\pi}^{\frac12\pi} e^{i\omega t}\ \psi(t)\,dt,\quad \psi(t)=e^{-1/\cos t},
\end{equation}
can be expanded as in \eqref{eq:mstph05}, and all terms will vanish. We conclude that $G(\omega)$ is exponentially small, but it is not clear if we have a simple estimate $\bigO\left(e^{-\omega}\right)$.

Another example is
\begin{equation}\label{eq:mahl1}
I(a,b)=\int_a^b\left(\frac{x+a}{x-a}\right)^{2ai}
\left(\frac{b-x}{b+x}\right)^{2bi}\,\frac{dx}x;
\end{equation}
$a$ and $b$ are large positive parameters with $b/a=c$, a 
constant
greater than 1. $I(a,b)$ is used in \cite[Eq.~(6.64)]{Lekner:1987:TRE} for describing the
Rayleigh approximation for a reflection amplitude in the theory
of electromagnetic and particle waves.

Because of the nature of the singular endpoints this integral is of a different type considered so far, but we may try to use the method of stationary phase. This is done in \cite{Mahler:1986:OSI}, where it was shown that the integral is exponentially small.  As Mahler
admitted, his results do not imply estimates of the form 
$I(a,b)=\bigO(e^{-a})$.
In this example the integrand can be considered for complex values of $x$, and by modifying the contour, leaving the real interval and entering the complex plane, we can show \cite{Temme:1989:AES} that  $I(a,b)=\bigO(e^{-2\pi a})$.

\begin{remark}\label{rem:trouble}
In fact, this is one of the troublesome aspects of the method of stationary phase (the use of neutralizers is another issue). The integral in \eqref{eq:mstph01} is defined on the real interval $[a,b]$ with functions $\phi,\psi\in C^n[a,b]$, possibly with $n=\infty$, and $\phi$ and $\omega$ real. When the integral is exponentially small, as $G(\omega)$ given in \eqref{eq:mstph37}, the method cannot give more information, unless we can extend the functions $\phi$ and $\psi$ to analytic functions defined in part of the complex plane. In that case more powerful methods of complex analysis become available, including modifying the interval into a contour  and using the saddle point method.
 \eoremark
 \end{remark}

\section{Uniform expansions}\label{sec:unimeth}   

The integrals considered so far are in the form
\begin{equation}\label{eq:stfm01}
\int_\calC e^{-\omega\phi(z)}\psi(z)\,dz
\end{equation}
along a real or complex path $\calC$ with one large parameter $\omega$. In many problems in physics, engineering, probability theory, and so on, additional parameters are present in $\phi$ or $\psi$, and these parameters may have influence on the validity of the asymptotic analysis, the rate of asymptotic convergence of the expansions, or on  choosing  a suitable method or contour. For example, complications will arise when $\phi$ or $\psi$ have
singular points that move to a saddle point under the influence of extra parameters. 
But it is also possible that two decisive points, say, two saddle points are proximate or coalescing. 

We will explain some aspects of uniform asymptotic expansions, usually with cases that are relevant in the asymptotic behavior of special functions.

\subsection{Van der Waerden's method}\label{sec:waerden}
In 1951 B.~L.~van der Waerden\footnote{Van der Waerden was appointed professor in the University of Amsterdam  from 1948--1951, and stayed also at the Mathematisch Centrum during that period.} wrote a paper \cite{Waerden:1951:OMS} in which he demonstrated what to do for integrals of the type (his notation was different)
\begin{equation}\label{eq:waerden01}
F_\alpha(\omega)=\frac{1}{2\pi i}\intr e^{-\omega t^2}\frac{f(t)}{t-
i\alpha}\,dt, \quad \Re\omega>0,
\end{equation}
where we assume that $f$ is analytic in a strip $\calD_a=\{t\in\CC \colon \vert\Im t\vert \le a, \Re t\in\RR\}$ for some positive $a$. Initially we take $ \Re\alpha >0$ but to let the pole cross the real axis we can modify the contour to avoid the pole, or use residue calculus.

Because $\vert\alpha\vert$ may be small we should not expand $f(t)/(t-i\alpha)$ in powers of $t$ as in Laplace's method (see \S\ref{sec:laplacemethod}). To describe Van der Waerden's method we write 
\begin{equation}\label{eq:waerden05}
f(t)=\left(f(t)-f(i\alpha)\right)+f(i\alpha).
\end{equation}
Then the integral in \eqref{eq:waerden01} can be written in the form
\begin{equation}\label{eq:waerden06}
F_\alpha(\omega)=\tfrac12f(i\alpha) e^{\omega\alpha^2}\erfc\left(\alpha\sqrt{{\omega}}\right)+
\frac{1}{2\pi i}\intr e^{-\omega t^2}g(t)\,dt,
\end{equation}
where
\begin{equation}\label{eq:waerden07}
g(t)=\frac{f(t)-f(i\alpha)}{t-i\alpha}.
\end{equation}
We have used one of the  error functions
\begin{equation}\label{eq:waerden02}
\erf\,z =\frac{2}{\sqrt{\pi}}\int_0^z e^{-t^2}\,dt,\quad
\erfc\,z =\frac{2}{\sqrt{\pi}}\int_z^\infty e^{-t^2}\,dt,
\end{equation}
with the properties 
\begin{equation}\label{eq:waerden03}
\erf\,z+\erfc\,z=1,\quad \erf(-z)=-\erf\,z,\quad \erfc(-z)=2-\erfc\,z.
\end{equation}
These functions are in fact the normal distribution functions from probability theory.
But we also have the representation
\begin{equation}\label{eq:waerden04}
\frac{e^{-\omega\alpha^2}}{2\pi i}\intr e^{-\omega
t^2}\frac{dt}{t-i\alpha}=\tfrac12\erfc\left(\alpha\sqrt{{\omega}}\right),\quad \Re\alpha >0,
\end{equation}
and for other values of $\alpha$ we can use analytic continuation. To prove this representation, differentiate the left-hand side with respect to 
$\omega$.

The functions $f$ and $g$ of \eqref{eq:waerden07} are analytic  in the same domain $\calD_a$ and $g$ can be expanded in powers of $t$. When we substitute $\dsp{g(t)=\sn c_n(\alpha)t^n}$ into \eqref{eq:waerden06}, we obtain the large $\omega$ asymptotic representation of $F_\alpha(\omega)$ by using Laplace's method. This gives
\begin{equation}\label{eq:waerden08}
F_\alpha(\omega)\sim\tfrac12f(i\alpha) e^{\omega\alpha^2}\erfc\left(\alpha\sqrt{{\omega}}\right)+
\frac{1}{2 i\sqrt{\pi\omega}}\sn \left(\tfrac12\right)_n \frac{c_{2n}(\alpha)}{\omega^n}.
\end{equation}
A special feature is that this expansion also holds when $\alpha=0$ and when $\Re\alpha<0$. All coefficients $c_n(\alpha)$ are well-defined and analytic for  $i\alpha\in\calD_a$. In fact the expansion in \eqref{eq:waerden08} holds uniformly with respect to $i\alpha$ in a domain $\calD^*$ properly inside $\calD_a$. 

Van der Waerden used this method in a problem of Sommerfeld concerning the propagation of radio waves over a plane earth. He obtained a simpler uniform expansion than the one given by Ott \cite{Ott:1943:DSU}, who obtained a uniform expansion in which each term is an incomplete gamma function. 

\subsubsection{An example from De Bruijn's book}\label{sec:bruijn}
In De Bruijn's book \cite[\S5.12]{Bruijn:1981:AMA}  the influence of poles near the saddle point  is considered by studying the example
\begin{equation}\label{eq:bruijn01}
F_\alpha(z)=\beta^2\intr e^{-\omega t^2}\frac{f(t)}{\beta^2+ t^2}\,dt,\quad \beta=\omega^{-\frac12\alpha}, \quad\alpha>0,
\end{equation}
where $\omega$ is a positive large parameter. Observe that for all $\alpha$ the parameter $\beta $ is small.

Three separate cases are distinguished: $0<\alpha<1$, $\alpha=1$, and $\alpha>1$, for the special choice $f(t)=e^t$. For each case an asymptotic expansion is given. These expansions are really different in the sense that they do not pass into each other when $\alpha$ passes unity. Splitting off the pole is not considered.

We can use partial fraction decomposition to get two integrals with a single pole, but we can also write
\begin{equation}\label{eq:bruijn02}
f(t)=a_0+b_0t+\left(\beta^2+t^2\right)g(t),
\end{equation}
where we assume that $g$ is regular at the points $\pm i\beta$. This gives 
\begin{equation}\label{eq:bruijn03}
a_0=\frac{f(i\beta)+f(-i\beta)}{2},\quad
b_0=\frac{f(i\beta)-f(-i\beta)}{2i},
\end{equation}
and where $g(t)$ follows with these values from \eqref{eq:bruijn02}. 

Hence,
\begin{equation}\label{eq:bruijn04}
F_\alpha(z)=a_0\beta^2\intr e^{-\omega t^2}\frac{dt}{\beta^2+ t^2}\,dt+\beta^2\intr e^{-\omega t^2} g(t)\,dt,
\end{equation}
where the first integral can be written in terms of the complementary error function (see \eqref{eq:waerden04}):
\begin{equation}\label{eq:bruijn05}
F_\alpha(z)=a_0\beta\pi \,e^{\omega^{1-\alpha}}\erfc\left(\omega^{\frac12(1-\alpha)}\right)+\beta^2\intr e^{-\omega t^2}g(t)\,dt.
\end{equation}
After expanding $\dsp{g(t)=\sk g_k(\beta)t^k}$ the asymptotic representation follows:
\begin{equation}\label{eq:bruijn06}
F_\alpha(z)\sim a_0\beta\pi \,e^{\omega^{1-\alpha}}\erfc\left(\omega^{\frac12(1-\alpha)}\right)+
\beta^2\sqrt{\frac{\pi}{\omega}}\sk c_k(\beta)\frac{1}{\omega^k},
\end{equation}
where $\beta$ is defined in \eqref{eq:bruijn01} and
\begin{equation}\label{eq:bruijn07}
c_k(\beta)= g_{2k}(\beta)\left(\tfrac12\right)_k.
\end{equation}
These and all other coefficients are regular when $\beta\to0$. The expansion in \eqref{eq:bruijn06} is valid for all $\alpha\ge0$.

We give a few coefficients for $f(t)=e^t$, which is used in \cite{Bruijn:1981:AMA}. We have
\begin{equation}\label{eq:bruijn08}
\begin{array}{ll}
\dsp{c_0(\beta)= \frac{1-\cos\beta}{\beta^2}},\\[8pt]
\dsp{c_1(\beta)= \frac{\beta^2-2+2\cos\beta}{4\beta^4}},\\[8pt]
\dsp{c_2(\beta)= \frac{\beta^4-12\beta^2+24-24\cos\beta}{32\beta^6}}.
\end{array}
\end{equation}

From the representation in \eqref{eq:bruijn05} we see at once the special value $\alpha=1$: when $0<\alpha<1$ the complementary error function can be expanded by using the asymptotic expansion
\begin{equation}\label{eq:bruijn09}
\erfc\,z\sim\frac{e^{-z^2}}{z\sqrt{\pi}}\left(1-\frac{1}{2z^2}+\frac{3}{4z^4}-\frac{15}{8z^6}+\cdots\right),\quad z \to \infty.
\end{equation}
When $\alpha>1$ it can be expanded by using the the convergent power series expansion of $\erf\,z=1-\erfc\,z$.
When $\alpha=1$, De Bruijn gives an expansion in terms of functions related with the complementary error function. His first term is $\beta\,\pi\, e\,\erfc(1)$, which corresponds to our term $\beta\cos\beta\,\pi\, e\,\erfc(1)$.

It should be noted that De Bruijn was not aiming to obtain a uniform expansion with respect to $\alpha$. His discussion was about the  \lq\lq{range of the saddle point}\rq\rq, that is, the role of an extra parameter which causes poles in the neighborhood of the saddle point. But it is remarkable that he has not referred to Van der Waerden's method, which gives a very short treatment of the lengthy discussions that De Bruijn needs in his explanations.  

\begin{remark}\label{rem:bruijn01}
The coefficients in \eqref{eq:bruijn08} are well defined as $\beta\to0$. We see that, and this often happens  in uniform expansions, close to some special value of a parameter (in this case $\beta=0$), the numerical evaluation of the coefficients is not straightforward, although they are defined properly in an analytic sense. In the present case, with the special example $f(t)=e^t$ in \eqref{eq:bruijn01}, it is rather easy to compute the coefficients by using power series expansions in $\beta$, but for more general cases special numerical methods are needed.
\eoremark
\end{remark}

\subsection{The Airy-type expansion}\label{sec:airy}

In 1957, Chester, Friedman and Ursell published a pioneering paper \cite{Chester:1957:ESD} in which Airy functions were used as leading terms in expansions. They extended the method of saddle points for the case that two saddle points, both relevant for describing the asymptotic behavior of the integral, are nearby or even coalescing.

\begin{figure}[tb]
\vspace*{0.8cm}
\begin{center}
\begin{minipage}{5.5cm}
\epsfxsize=5.5cm \epsfbox{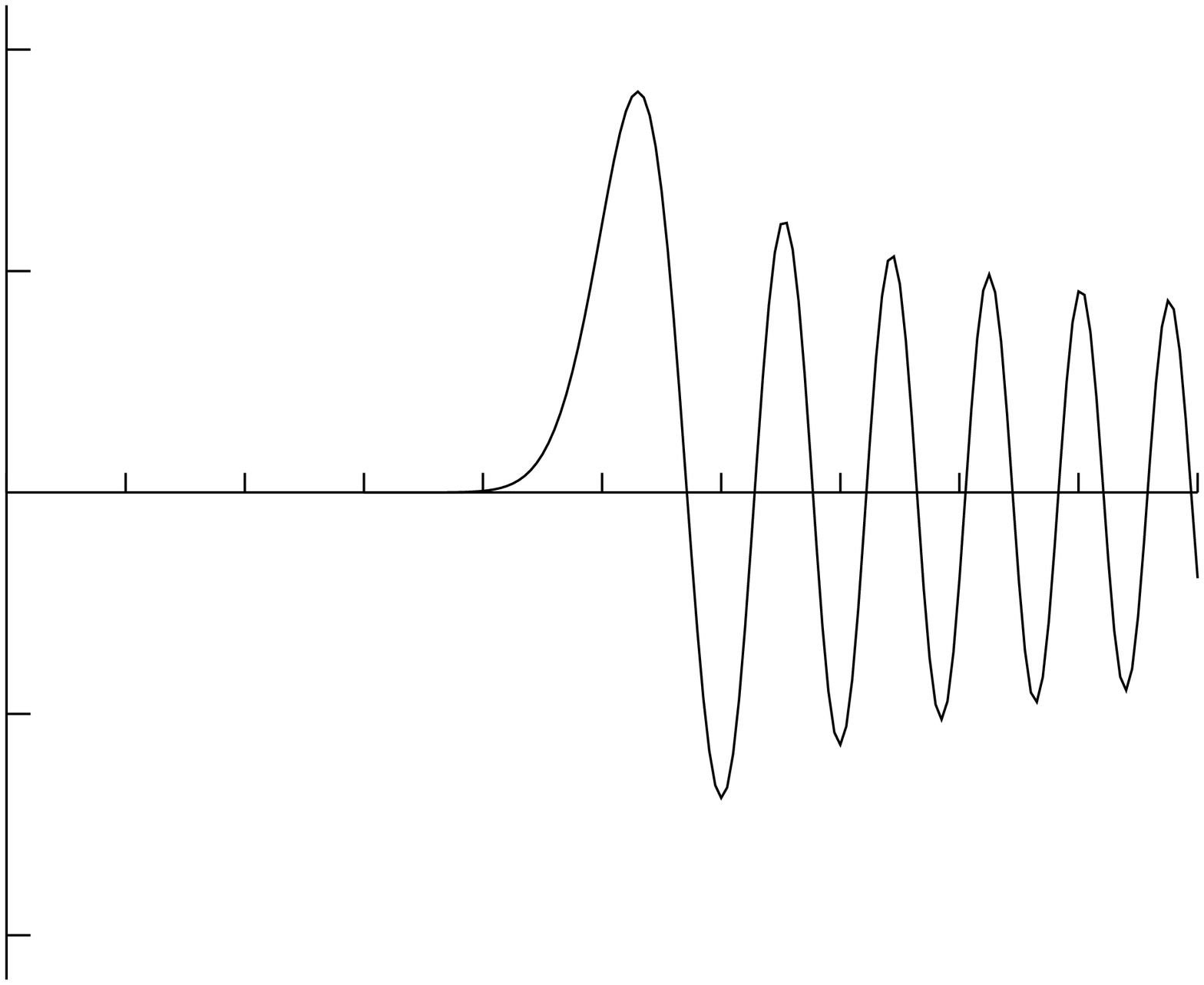}
\end{minipage}
\hspace*{0.5cm}
\begin{minipage}{5.5cm}
\epsfxsize=5.5cm \epsfbox{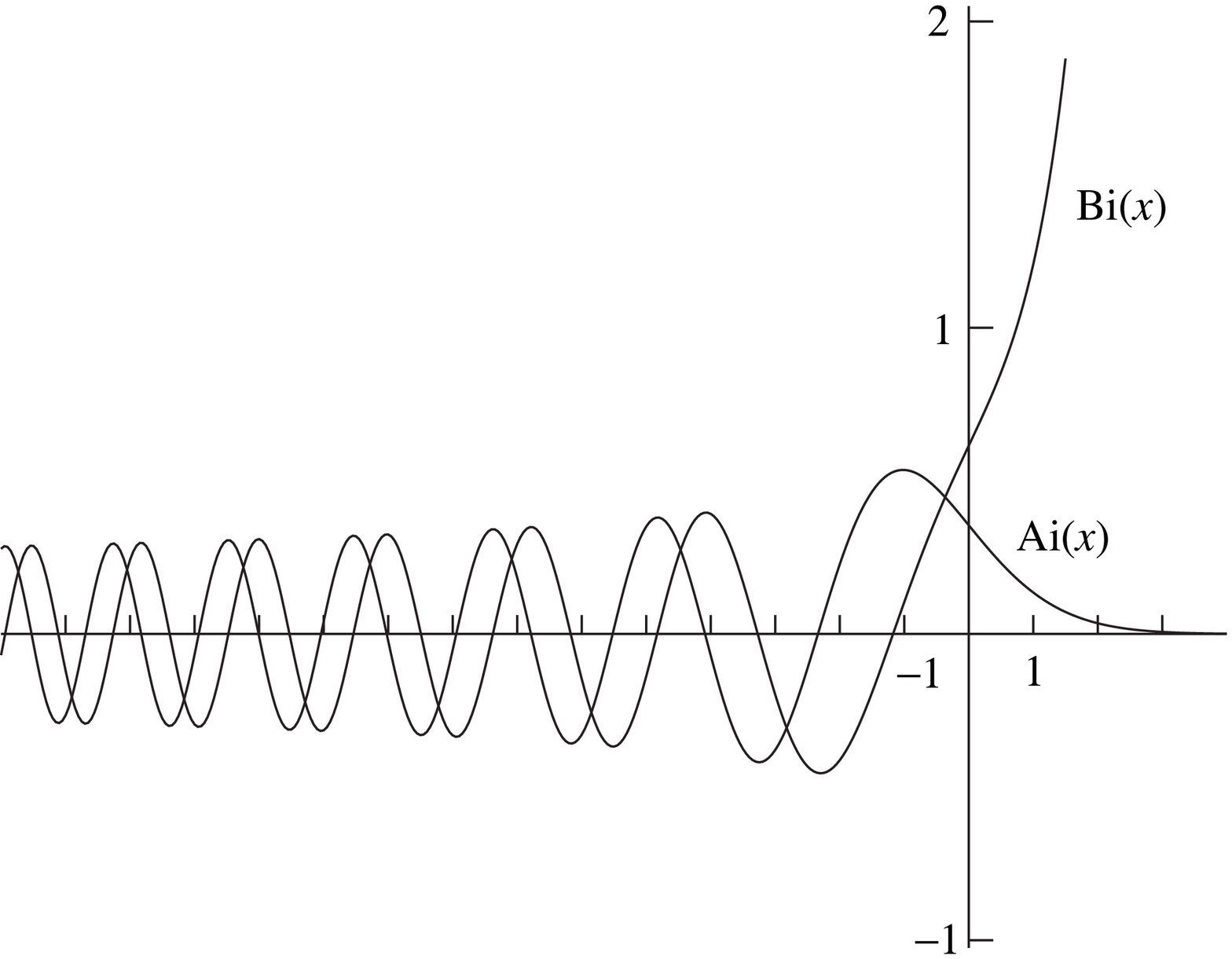}
\end{minipage}
\end{center}
\caption{{\bf Left:} The Bessel function $J_{50}(x), 0\le x\le100$. {\bf Right:} Graphs of the Airy functions $\Ai(x)$ and $\Bi(x)$ on the real line.
}
\label{fig:AiryBess}
\end{figure}

We describe this phenomenon for an integral representation of the Bessel function $J_\nu(z)$. We have
\begin{equation}\label{eq:airy02}
J_\nu(\nu z)=\frac{1}{2\pi i}\int_\calC
e^{\nu\phi(s)}\,ds, \quad \phi(s)=z\sinh s-s,
\end{equation}
where the contour $\calC$ starts at $\infty-\pi i$ and terminates at $\infty+\pi i$. We consider positive values of $z$ and $\nu$, with $\nu$ large.
From graphs of the 
Bessel
function of high order, see Figure~\ref{fig:AiryBess},  it can be seen that $J_\nu(x)$ starts
oscillating when $x$ crosses the value $\nu$. We concentrate on 
the transition area $x\sim \nu$. The Airy function $\Ai(-x)$ shows a similar behavior, as also follows from Figure~\ref{fig:AiryBess}.

There are two saddle points defined by the equation $\cosh s=1/z$. For $0<z\le1$ the saddle points are real, and  in that case we have
\begin{equation}\label{eq:airy03}
s_\pm= \pm\theta,\quad {\rm where}\quad z=1/\cosh \theta,\quad \theta\ge0.
\end{equation}
 When $0<z\le z_0<1$, with $z_0$ a fixed number, the positive saddle point $s_+$ is the relevant one. We can use the transformation $\phi(s)-\phi(s_+)=-\frac12t^2$ and apply Laplace's method of \S\ref{sec:laplacemethod}. This gives Debye's expansion for $J_\nu(z)$, see \cite[\S10.19(ii)]{Olver:2010:BFS}.

If $z$ is close to 1, both saddle point are relevant. As shown for the first time in \cite{Chester:1957:ESD}, we can use a cubic polynomial by writing
\begin{equation}\label{eq:airy07}
\phi(s)=\tfrac13 t^3-\zeta t+A,
\end{equation}
where $\zeta$ and $A$ follow from the condition that the saddle points  $s_\pm$ given in \eqref{eq:airy03}  correspond to $\pm\sqrt\zeta$ (the saddle points in the $t-$plane). It is not difficult to verify that $A=0$ (take $s=t=0$) and that
\begin{equation}\label{eq:airy08}
\tfrac23\zeta^{\frac32}=\ln\frac{1+\sqrt{1-z^2}}{z}-\sqrt{1-z^2},\quad 0<z\le 1.
\end{equation}
This relation can be analytically continued for $z>1$ and for complex values of $z$, $z\ne-1$. For small values of $\vert \zeta\vert$ we have 
\begin{equation}\label{eq:airy09}
z(\zeta)= \sn z_n\eta^n=1-\eta+\tfrac3{10}\eta^2+\tfrac1{350}\eta^3- 
\tfrac{479}{63000}\eta^4+\ldots,\quad \eta= 2^{-\frac13}\zeta,
\end{equation}
with convergence for $\vert\zeta\vert < (3\pi/2)^{2/3}$.

By using conditions on the mapping in \eqref{eq:airy07} it can be proved that \eqref{eq:airy02} can be written as
\begin{equation}\label{eq:airy10}
J_\nu(\nu z)=\frac{1}{2\pi i}\int_{\infty e^{-\frac13\pi i}}^{\infty e^{\frac13\pi i}}e^{\nu(\frac13t^3-\zeta 
t)}g(t)\,dt,
\quad g(t)=\frac{ds}{dt}=\frac{t^2-\zeta}{\phi^\prime(s)}.
\end{equation}
It is not difficult to verify that $g$ can be defined at both saddle points $t=\pm\sqrt\zeta$. We have
\begin{equation}\label{eq:airy11}
g\left(\pm\sqrt\zeta\right)=\left(\frac{4\zeta}{1-z^2}\right)^{\frac{1}{4}}.
\end{equation}

A detailed analysis is needed to show that  the transformation in \eqref{eq:airy07} one-to-one maps  a relevant part of the $s-$plane to the  $t-$plane. One of the three solutions of the cubic equation has to be chosen, and the proper one should be such that we can write the original integral in \eqref{eq:airy02}  along a chosen contour in the form of \eqref{eq:airy10}, with the contour at infinity as indicated. The proofs in  \cite{Chester:1957:ESD} are quite complicated and later more accessible methods came available, in particular when differential equations are used as starting point of the asymptotic analysis. For details of the form of the expansion we refer to \cite[\S10.20]{Olver:2010:BFS}, where also references for proofs are given. 

In \cite[\S VII.5]{Wong:2001:AAI} a detailed explanation is given of the various aspects of the conformal mapping for such transformatio for the case of an Airy-type expansion of the Laguerre polynomials starting with an integral representation. 

Using the construction explained  \cite[\S VII.4]{Wong:2001:AAI} we can find for the integral in \eqref{eq:airy10}
the representation
\begin{equation}\label{eq:twocoal08}
J_\nu(\nu z)= g\left(\sqrt{\zeta}\right) \left(
\frac{\Ai\left(\zeta\nu^{\frac23}\right)}{\nu^{\frac13}}A(\nu,\zeta)+
 \frac{\Ai^{\prime}\left(\zeta\nu^{\frac23}\right)}{\nu^{\frac53}}B(\nu,\zeta)\right),
\end{equation}
where $\zeta$ is defined in \eqref{eq:airy08} and  $g\left(\sqrt{\zeta}\right)$ is given in \eqref{eq:airy11}, and
\begin{equation}\label{eq:twocoal09}
A(\nu,\zeta)\sim\sn\frac{A_n(\zeta)}{\nu^{2n}},\quad B(\nu,\zeta)\sim \sn \frac{B_n(\zeta)}{\nu^{2n}}.
\end{equation}
The asymptotic expansions are valid for large $\nu$ and uniformly for $z\in(0,\infty)$. 
See \cite[\S10.20(i)]{Olver:2010:BFS} for more details on the coefficients.

\subsubsection{A few remarks on Airy-type expansions}\label{sec:airyrem}
Airy functions are special cases of Bessel functions of order $\pm\frac13$ and are  
named after G.B.~Airy (1838), a British astronomer, who used them in the study of rainbow phenomena. They occur in many other problems from physics, for example as solutions to boundary value problems in quantum mechanics and electromagnetic theory. 

The Airy function $\Ai(z)$ can be defined by
\begin{equation}\label{eq:twocoal02}
\Ai(z)=
\frac1{2\pi i}\int_{\infty e^{-\frac13\pi i}}^{\infty e^{\frac13\pi i}}e^{\frac13t^3-z t}\,dt,
\end{equation}
and we see that the integral in \eqref{eq:airy10} becomes an Airy function when the function $g$ is replaced with a constant.

Airy functions are solutions of Airy's differential equation $\dsp{\frac{d^2w}{dz^2}=zw}$.
It is not difficult to verify that the function in \eqref{eq:twocoal02} satisfies this equation. Two independent solutions are denoted by $\Ai(z)$ and $\Bi(z)$, which are entire functions and real when $z$ is real.  
They  are oscillatory  for $z<0$ and decrease ($\Ai$) or increase ($\Bi$) exponentially fast for $z>0$; see Figure~\ref{fig:AiryBess}.
 
Airy's equation is the simplest second-order linear differential equation
showing a turning point (at $z=0$). A turning point in a differential equation (see  \cite[Chapter~11]{Olver:1997:ASF}) corresponds to  integrals that have two coalescing saddle points, as we have seen for the Bessel function.

Many special functions have Airy-type expansions. For example, the orthogonal polynomials oscillate inside the domain of orthogonality and become non-oscillating outside this domain. The large degree behavior of these polynomials can often be described by Airy functions, in particular to describe the change of behavior near finite endpoints of the orthogonality domain. 

In the example of the Bessel function $J_\nu(z)$, a function of two variables,  we use the Airy functions $\Ai(z)$, a function of one variable. This is a general principle when constructing uniform expansions in which special functions are used as approximations: reduce the number of parameters.

The results on Airy-type expansions for integrals of Chester {\em et all.}~in 1957 were preceded by fundamental results  in 1954 of Olver, who used differential equations for obtaining these expansions, see \cite{Olver:1954:AEB,Olver:1954:ASL}. Olver extended earlier results for turning point equations and in his work realistic and computable error bounds are constructed for the remainders in the expansions.

The form given in \eqref{eq:twocoal08}  follows from a method of Bleistein \cite{Bleistein:1966:UAE}, who used it for a different type of integral. Chester {\em et all.}~have not given this form; they expanded the function $g$ in \eqref{eq:airy10} in a two-point Taylor series (see \cite{Lopez:2002:TPT}) of the form
\begin{equation}\label{eq:remairy01}
g(t)=\sn a_n\left(t^2-\zeta\right)^n+t\sn b_n\left(t^2-\zeta\right)^n,
\end{equation}
and the expansion obtained in this way can be rearranged to get the expansion given in \eqref{eq:twocoal08}.

The paper by Chester {\em et all.}~\cite{Chester:1957:ESD} appeared in 1957 (one year before De Bruijn published his book), and together with Van der Waerden's paper \cite{Waerden:1951:OMS} it was the start of 
a new period of asymptotic methods for integrals. As mentioned earlier, Airy-type expansions were considered in detail by Olver (see also \cite[Chapter~11]{Olver:1997:ASF}), who provided realistic error bounds for remainders by using differential equations. In \cite{Olde:1994:UAE} it was shown how to derive estimates of remainders in Airy-type expansions obtained from integrals.

It took some time to see how to use  three-term recurrence relations or linear  second order difference equations as starting point for obtaining Airy-type expansions or other uniform expansions. The real breakthrough came in a paper by Wang and Wong \cite{Wang:2002:UAE}, with an application for the Bessel function $J_\nu(z)$. But the method can also be applied to cases when no differential equation or integral representation is available, such as  for polynomials associated with the Freud weight $e^{-x^4}$ on $\RR$; see \cite{Wang:2003:AEF}. In \cite{Deift:1993:ASD} a completely new method was introduced. This
method is based on Riemann-Hilbert type of techniques. All the methods mentioned
above can be applied to obtain Airy-type expansions but also to a number of other cases. In the next section we give more examples of such standard forms.

\section{A table of standard forms}\label{sec:standard}

\renewcommand{\arraystretch}{2.5}
  \begin{table}
  \caption{An overview of standard forms\label{tab:overstand}}
 \begin{center}
\begin{tabular}{||r|c|c|c||}
\hline
\bf No. &\bf Standard Form& \bf Approximant& \bf Decisive points \\
\hline
\rm
1&$\dsp{\int_{0}^{\infty}e^{-zt}\frac{f(t)}{t+\alpha}\,dt}$&Exponential integral 
&$0,-\alpha$\\
2&$\dsp{\int_{-\infty}^{\infty}e^{-zt^2}\frac{f(t)}{t-i\alpha}\,dt}$&Error 
function&$0,i\alpha$\\
3&$\dsp{\int_{-\infty}^{\alpha}e^{-
zt^2}f(t)\,dt}$&Error 
function&$0,\alpha$\\ 
4&$\dsp{\int_0^\infty t^{\beta-1}e^{-z(\frac12t^2-\alpha 
t)}f(t)\,dt}$& Par.~cylinder
function&$0,\alpha$\\
5&$\dsp{\int_{-\infty}^{\infty}e^{-zt^2}\frac{f(t)}{(t-
i\alpha)^\mu}\,dt}$& Par.~cylinder
function&$0,i\alpha$\\
6&$\dsp{\int_{\calL}e^{z(\frac13t^3-\alpha t)}f(t)\,dt}$& Airy 
function&$\pm\sqrt{{\alpha}}$\\
7&$\dsp{\int_0^\infty t^{\lambda-1}e^{-zt}f(t)\,dt}$ &Gamma 
function&$0,\lambda/z$\\
8&$\dsp{\int_\alpha^\infty t^{\beta-1}e^{-zt}f(t)\,dt}$ &Inc. gamma
function&$0,\alpha,\beta/z$\\
9&$\dsp{\int_{-\infty}^{(0+)} t^{-\beta-1}e^{z(t+\alpha/t)}f(t)\,dt}$& Bessel $I$
function&$0,\pm\sqrt{{\alpha}}$\\
10&$\dsp{\int_0^\infty t^{\beta-1}e^{-z(t+\alpha/t)}f(t)\,dt}$& Bessel $K$
function&$0,\pm\sqrt{{\alpha}}$\\
11&$\dsp{\int_0^\infty t^{\lambda-1}(t+\alpha)^{-\mu}e^{-zt}f(t)\,dt}$& Kummer $U$ 
function&$0,-\alpha$\\
12&$\dsp{\int_{-\alpha}^\alpha e^{-zt}(\alpha^2-t^2)^\mu f(t)\,dt}$& Bessel $I$
function&$0,\pm\alpha$\\ 
13&$\dsp{\int_\alpha^\infty e^{-zt}(t^2-\alpha^2)^\mu f(t)\,dt}$& Bessel $K$
function&$0,\pm\alpha$\\ 
14&$\dsp{\int_0^\infty\frac{\sin z(t-\alpha)}{t-\alpha}f(t)\,dt}$& Sine
integral&$0,\alpha$\\[6pt]
\hline
\end{tabular}
  \end{center}
  \end{table}
\renewcommand{\arraystretch}{1.0}

In Table~\ref{tab:overstand} we give an overview of standard forms considered in the literature. In almost all cases of the table, these forms arise in the asymptotic analysis of special functions,  and in all cases special functions are used as leading terms in the approximations. 

The decisive points mentioned in the table are the points in the interval of integration (or close to this interval) where the main contributions to the approximation can be obtained. When these points are coalescing, uniform methods are needed.

The cases given in the table can be used to solve asymptotic problems for special functions, for problems in science and engineering, in probability theory, in number theory, and so on. We almost never see the forms as shown in the table ready-made in these problems, but we need insight, technical steps, and nontrivial transformations to put the original integrals, sums, solutions of equations,  and so on, into one of these standard forms.

We give a few comments on the cases in Table~\ref{tab:overstand}.
\begin{description}
\item[Case 1]
We can use Van der Waerden's method, see \S\ref{sec:waerden}. When $f=1$ the integral reduces to the exponential integral. This simple case has not been considered in the literature and we give an example. Let $\dsp{S_n(z)=\sum_{k=1}^n \frac{z^k}{k}}$. Determine the behavior as $n\to\infty$,
that holds uniformly for values of $z$ close to 1. 
Observe that as $n\to\infty$,
\begin{equation}\label{eq:cases01}
S_n(1)\sim \ln n, \quad S_n(z)\sim-\ln(1-z), \quad  |z|<1,
\end{equation}
and that $S_n(z)$ has the representation
\begin{equation}\label{eq:cases02}
S_n(z)=-\ln(1-z)-z^{n+1}\intp\frac{e^{-ns}}{e^s-z}\,ds, \quad z\ne 1.
\end{equation}
There is a pole at $s=\ln z$, at the negative axis when $0<z<1$, which can be split off, giving the exponential integral, and this special function will deal with the term $-\ln(1-z)$ in \eqref{eq:cases02} as $z\to1$.

\item[Case 2]
This has been considered in \S\ref{sec:waerden}. 

\item[Case 3]
Again we can use the complementary error function. This standard form has important applications for cumulative distribution functions of probability theory. As explained in \cite{Temme:1982:UAE}, the well-known gamma and beta distributions, and several other ones, can be transformed into 
this standard form.

\item[Case 4]
When $f=1$, the integral becomes the parabolic cylinder function $U(a,z)$  and when  $\beta=1$  this case is the same as Case 3. See  \cite{Bleistein:1966:UAE}, where an integration-by-parts method was introduced that can be used in slightly different versions to obtain many other uniform expansions. See also \cite[Chapter~7]{Wong:2001:AAI}.

\item[Case 5]
Again we use a parabolic cylinder function with $\alpha$ in a domain containing the origin. If needed the path of integration should be taken around the branch cut. For an application to Kummer functions, see \cite{Temme:1978:UAE}.

\item[Case 6]
This case has been considered in \S\ref{sec:airy}.

\item[Case 7]
We assume that $z$ is large and that $\lambda$ may be large as well. This is different from Watson's lemma considered in \S\ref{sec:aslap}, where we have assumed that $\lambda$ is fixed. There is a saddle point at $\mu=\lambda/z$.  For details we refer to \cite{Temme:1983:LAP,Temme:1985:LTI}.

\item[Case 8]
This case extends the previous one with an extra parameter $\alpha$, and we assume
$\lambda\ge0$, $\alpha\ge0$ and $z$ large. As in Case~7,  $\lambda$ may also 
be large; $\alpha$ may be large as well, even larger than $\lambda$. When $f=1$ this integral becomes an incomplete gamma function. 
In \cite{Temme:1987:ILI} we have given more details with an application to the incomplete beta functions.

\item[Case 9]
When $f=1$ the contour integral reduces to a modified $I-$Bessel function when $\alpha>0$, and to a $J-$Bessel function when $\alpha<0$. For an application to Laguerre polynomials, we refer to \cite[Chapter~7]{Wong:2001:AAI}.

\item[Case 10]
This  real integral  with the essential singularity at  the origin is related with the previous case, and
reduces to a modified $K-$Bessel function when $f$ is a constant. We can give an asymptotic expansion of this integral for large values of $z$, which is uniformly valid for $\lambda\ge0$, $\alpha\ge0$. For more details we refer to \cite{Temme:1990:UAE}, where an application is given to confluent hypergeometric functions.

\item[Case 11]
This is more general  than Case 1. See \cite{Oberhettinger:1959:OAM} for a more details.

\item[Case 12 and Case 13]
When $f(t)=1$  the integrals reduce to modified Bessel functions, see
\cite{Ursell:1984:IWL} (with applications to Legendre functions). In 
\cite{Ursell:2007:IWN} a contour integral is considered, with an application to Gegenbauer polynomials. In that case the expansion is in terms of the $J-$Bessel function.

\item[Case 14]
This is considered in \cite{Zilbergleit:1977:UAE}, where a complete asymptotic expansion is derived in which the sine integral is used for a smooth transition of $x=0$ to $x>0$.

\end{description}

\section{Concluding remarks}\label{sec:concl}
We have given an overview of the classical methods for obtaining asymptotic expansions for integrals, occasionally by giving comments on De Bruijn's book \cite{Bruijn:1981:AMA}. We have explained how uniform expansions for integrals were introduced by Van der Waerden \cite{Waerden:1951:OMS} and Chester {\em et all.} \cite{Chester:1957:ESD}, and we have given an overview of uniform methods that were introduced since the latter pioneering paper. For Airy-type expansions we have mentioned several other approaches, which are also  used for other expansions in which special functions are used as leading terms in the approximations.

\section*{Acknowledgements}
The author thanks the referee for advice and comments; he acknowledges research facilities from CWI, Amsterdam, and support from  {\emph{Ministerio de Ciencia e Innovaci\'on}, Spain}, 
project MTM2009-11686.

\def\cprime{$'$} \def\cprime{$'$}

{

}
\end{document}